\documentclass[12pt]{article}

\usepackage{amsmath}
\usepackage{amssymb}
\usepackage{graphicx}
\usepackage{epsfig}

\usepackage[margin=2cm]{geometry}
\newtheorem{theorem}{Theorem}
\newtheorem{lemma}[theorem]{Lemma}

\def \Z {{\mathbb{Z}}}

\def \a {\alpha}
\def \b {\beta}

\def \o {\omega}

\def \op {{\rm OP}}
\def \proof {\noindent{\bf Proof}\quad}

\def \P {\mathcal P}

\def \md#1{{\,({\rm mod}\ #1)}}

\def\endmark{\hskip 2em$\square$\par}
\def \qed {\hfill\endmark}

\def \cay {{\rm Cay}}

\title{\bf
On factorisations of complete graphs into
circulant graphs and the
Oberwolfach Problem}

\author{
Brian Alspach
\thanks{
School of Mathematical and Physical Sciences,
University of Newcastle,
Callaghan, NSW 2308, Australia.
\texttt{brian.alspach@newcastle.edu.au}
}
\and
Darryn Bryant
\thanks{
School of Mathematics and Physics,
The University of Queensland,
Qld 4072, Australia.
\texttt{db@maths.uq.edu.au},
}
\and
Daniel Horsley
\thanks{
School of Mathematical Sciences,
Monash University,
Vic 3800, Australia.
\texttt{danhorsley@gmail.com},
}
\and
Barbara Maenhaut
\thanks{
School of Mathematics and Physics,
The University of Queensland,
Qld 4072, Australia.
\texttt{bmm@maths.uq.edu.au}}
\and
Victor Scharaschkin
\thanks{
School of Mathematics and Physics,
The University of Queensland,
Qld 4072, Australia.
\texttt{victors@maths.uq.edu.au}
}
}

\date{ }

\begin{document}
\maketitle\thispagestyle{empty}
\def\baselinestretch{1.5}\small\normalsize

\begin{abstract}
Various results on factorisations of complete graphs into circulant graphs and on 2-factorisations of these circulant graphs are proved.
As a consequence, a number of new results on the Oberwolfach Problem are obtained. For example, a complete solution to the Oberwolfach Problem is given for every $2$-regular graph of order $2p$ where $p\equiv 5\md 8$ is prime.
\end{abstract}


\vspace{0.3cm}


\section{Introduction}

The Oberwolfach problem was posed by Ringel in the 1960s and is first mentioned in \cite{Guy}. It concerns graph factorisations. A {\em factor} of a graph is a spanning subgraph and a {\em factorisation} is a decomposition into edge-disjoint factors.
A factor that is regular of degree $k$ is called a {\em $k$-factor}. If each factor of a factorisation is a $k$-factor, then the factorisation is called a {\em $k$-factorisation}, and if each factor is isomorphic to a given graph $F$, then we say it is a {\em factorisation into $F$}.

Let $F$ be an arbitrary $2$-regular graph and let $n$ be the order of $F$.
If $n$ is odd, then the {\em Oberwolfach Problem} $\op(F)$ asks for a $2$-factorisation of $K_n$ into $F$, and if $n$ is even, then $\op(F)$ asks for a $2$-factorisation of $K_n-I$ into $F$, where $K_n-I$ denotes the graph obtained from $K_n$ by removing the edges of a $1$-factor.

The Oberwolfach Problem has been solved completely when $F$ consists of isomorphic components \cite{AlsHag,AlsSchStiWag,HofSch}, when $F$ has exactly two components \cite{Tra}, when $F$ is bipartite \cite{BryDan,Hag} and in numerous special cases. See \cite{BryRod} for a survey of results up to 2006. It is known that there is no solution to $\op(F)$ for
$F\in\{C_3\cup C_3, C_4\cup C_5, C_3\cup C_3\cup C_5, C_3\cup C_3\cup C_3\cup C_3\}$, but a solution exists for every other $2$-regular graph of order at most $40$ \cite{DezFraHuaMesRos}.

In \cite{BrySch}, it was shown that the Oberwolfach Problem has a solution for every $2$-regular graph of order $2p$ where $p$ is any of the infinitely many primes congruent to $5\md{24}$, and for every $2$-regular graph whose order is in an infinite family of primes congruent to $1\md{16}$. In this paper we extend these results as follows. We show that $\op(F)$ has a solution for every $2$-regular graph of order $2p$ where $p$ is any prime congruent to $5\md 8$ (see Theorem \ref{mainThmp=5mod8}), and we obtain solutions to $\op(F)$ for broad classes of $2$-regular graphs in many other cases (see Theorems \ref{combinedFactorsThm} and \ref{134result}). We also obtain results on the generalisation of the Oberwolfach Problem to factorisations of complete multigraphs into isomorphic $2$-factors (see Theorem \ref{lambdageq7thm}). Our results are obtained by constructing various factorisations of complete graphs into circulant graphs in Section \ref{SectionCompleteintoCircs}, and then showing in Section \ref{2FactorisationOfCirculantsSection} that these circulant graphs can themselves be factored into isomorphic $2$-regular graphs in a wide variety of cases.

\section{Factorising complete graphs into circulant graphs}\label{SectionCompleteintoCircs}

Let $G=(G,\cdot)$ be a finite group with identity $e$ and let $S$ be a subset of $G$ such that $e\notin S$ and $s\in S$ implies $s^{-1}\in S$. The {\em Cayley graph on $G$ with connection set $S$}, denoted $\cay(G\,;S)$, has the elements of $G$ as its vertices and $g$ is adjacent to $g\cdot s$ for each $s\in S$ and each $g\in G$. A Cayley graph on a cyclic group is called a {\em circulant graph}. We use the standard notation of $\Z_n$ for the ring of integers modulo $n$, and we use $\Z_n^*$ for the multiplicative group of units modulo $n$.

In this section we consider factorisations of $K_n$ for $n$ odd (in Section \ref{completeintocircs}) and of $K_n-I$ for $n$ even (in Section \ref{completeminusIintocircs}) into circulant graphs. A $2$-regular graph is a circulant if and only if its componenets are all isomorphic. Thus, for each $2$-regular circulant graph $F$, there exists a factorisation of $K_n$ (if $F$ has odd order) or of $K_n-I$ (if $F$ has even order) into $F$; except that there is no such factorisation when $F\in\{C_3\cup C_3,C_3\cup C_3\cup C_3\cup C_3\}$. Considerably less is known for factorisations into circulant graphs of degree greater than $2$. Some factorisations into $\cay(\Z_n\,;\pm\{1,2\})$ and $\cay(\Z_n\,;\pm\{1,2,3,4\})$ are given in \cite{Bry} and \cite{BrySch} respectively, and some further results, including results on self-complementary and almost self-complementary circulant graphs, appear in \cite{AlsMorVil,DobSaj,FroRosSir,PraLiStr}.

\subsection{Factorising complete graphs of odd order}\label{completeintocircs}

In this subsection we will construct factorisations of complete graphs of odd order into isomorphic circulant graphs by finding certain partitions of cyclic groups. Problems concerning such partitions have been well-studied, for example see \cite{Sza}, and existing results overlap with some of the results in this subsection. In particular, Theorem \ref{Kpinto123} below is a consequence of Lemma 3.1 of \cite{Mun}.

\begin{lemma}\label{intocircs}
Let $s$ be an integer,
let $p\equiv 1\md {2s}$ be prime, and let
$S=\pm\{d_1,d_2,\ldots,d_s\}\subseteq\Z_p^*$.
Further, suppose
$a$ and $b$ are integers such that $2abs=p-1$, let $G=(\Z_p^*)^b$, and let $H=(\Z_p^*)^{bs}$.
If
$d_1,d_2,\ldots,d_s$ represent the $s$ distinct cosets of $G/H$, then
there exists a $2s$-factorisation of $K_p$ into
$\cay(\Z_p\,;S)$.
\end{lemma}

\proof
For each $x\in\Z_p$ let
$xS=\{xy:y\in S\}$.
Since $p$ is prime, $\cay(\Z_p\,; xS)\cong
\cay(\Z_p\,;S)$ for any $x\in\Z_p\setminus\{0\}$.
If there is a partition of
$\Z_p^*$ into sets $x_1S,x_2S,\ldots,x_{ab}S$
where $x_i\in\Z_p\setminus\{0\}$ for
$i=1,2,\ldots,ab$, then
$\{\cay(\Z_p\,;x_iS):i=1,2,\ldots,ab\}$
is the required $2s$-factorisation of $K_p$.
We now present such a partition.

Let $\omega$ be a generator of $\Z_p^*$.
Thus, $H=\o^0,\o^{bs},\o^{2bs},\ldots,\o^{(2a-1)bs}$,
and $\o^{abs}=-1\in H$.
Let $A=\o^0,\o^{bs},\o^{2bs},\ldots,\o^{(a-1)bs}$,
so that $H=A\cup -A$ ($A$ is a set of representatives for the cosets in $H$ of the order $2$ subgroup of $H$).
Since $d_1,d_2,\ldots,d_s$ represent distinct cosets of $G/H$, it is easy to see that
$\{xS:x\in A\}$ is a partition of $G$.
Thus, if $B$ is a set of representatives
for the cosets of $\Z_p^*/G$,
then
$\{xyS:x\in A,y\in B\}$
is the required partition of $\Z_p^*$.
\qed

\vspace{0.3cm}

Note that upon putting $s=1$ in Lemma \ref{intocircs} we obtain the Hamilton decomposition $$\{\cay(\Z_p\,;\{\pm 1\}),\cay(\Z_p\,;\{\pm 2\}),\ldots,\cay(\Z_p\,;\{\pm\textstyle \frac{p-1}2\})\}$$ of $K_p$. We will be mostly interested in applications of Lemma \ref{intocircs} where the connection set $S$ is $\pm\{1,2\}$, $\pm\{1,2,3\}$, $\pm\{1,3,4\}$ or $\pm\{1,2,3,4\}$. The factorisations given by Lemma \ref{intocircs} have the property that each factor is invariant under the action of $\Z_p$. It is worth mentioning that for $S\in\{\pm\{1,2\},\pm\{1,2,3\},\pm\{1,3,4\},\pm\{1,2,3,4\}\}$, the construction given in Lemma \ref{intocircs} yields every $2s$-factorisation of $K_p$ into $\cay(\Z_p\,;S)$ with this property. This follows from the results in \cite{Bur} and \cite{Mom}, together with Turner's result \cite{Tur} that for $p$ prime $\cay(\Z_p\,;S)\cong\cay(\Z_p\,;S')$ if and only if there exists an $\a\in\Z_p^*$ such that $S'=\a S$.

\begin{theorem}
If $p\equiv 1\md 4$ is prime and $4$ divides the order of $k$ in $\Z_p^*$, then there is a factorisation of $K_p$ into $\cay(\Z_p\,;\pm\{1,k\})$.
\end{theorem}

\proof
Apply Lemma \ref{intocircs} with $S=\pm\{1,k\}$ taking $G$ to be the subgroup of $\Z_p^*$ generated by $k$, and $H$ to be the index $2$ subgroup of $G$.
\qed

\begin{theorem}\label{Kpinto123}
If $p\equiv 1\md{6}$ is prime such that $2,3\notin(\Z_p^*)^3$ and $6\in(\Z_p^*)^3$,
then there is a factorisation of $K_p$ into $\cay(\Z_p\,;\pm\{1,2,3\})$.
\end{theorem}

\proof
It follows from $2,3\notin(\Z_p^*)^3$ and $6\in(\Z_p^*)^3$ that $1$, $2$ and $3$ represent the three cosets of $\Z_p^*/(\Z_p^*)^3$. Thus, we obtain the required factorisation by applying Lemma \ref{intocircs} with $b=1$.
\qed

\begin{theorem}\label{Kpinto134}
If $p\equiv 1\md{6}$ is prime such that $2,3,6\notin(\Z_p^*)^3$,
then there is a factorisation of $K_p$ into $\cay(\Z_p\,;\pm\{1,3,4\})$.
\end{theorem}

\proof
It follows from $2,3,6\notin(\Z_p^*)^3$ that $1$, $3$ and $4$ represent the three cosets of $\Z_p^*/(\Z_p^*)^3$. Thus, we obtain the required factorisation by applying Lemma \ref{intocircs} with $b=1$.
\qed

\vspace{0.3cm}

The primes less than 1000 to which Theorem \ref{Kpinto123} applies are $$7,37,139,163,181,241,313,337,349,379,409,421,541,571,607,631,751,859,877,937,$$
and the primes less than 1000 to which Theorem \ref{Kpinto134} applies
are $$13,19,79,97,199,211,331,373,463,487,673,709,769,823,829,883,907.$$
In the next theorem we show that there are infinitely many primes to which Theorem \ref{Kpinto123} applies,
and also infinitely many primes to which Theorem \ref{Kpinto134} applies.

\begin{theorem}\label{Thm-p1mod6InfMany}
There are infinitely many values of $p$ such that $p$ is prime, $p\equiv 1\md 6$, $2,3\notin(\Z_p^*)^3$ and $6\in(\Z_p^*)^3$, and there are infinitely many values of $p$ such that $p$ is prime, $p\equiv 1\md 6$ and $2,3,6\notin(\Z_p^*)^3$.
\end{theorem}

\newcommand{\Q}{\mathbb{Q}}
\newcommand{\Zp}{\mathbb{Z}/p\mathbb{Z}}
\newcommand{\Fp}{\mathbb{F}_p}
\newcommand{\Fpx}{{\mathbb{F}_p^{\times}}}

\newcommand{\from}{\colon}
\renewcommand{\O}{\mathbb{O}} 
\newcommand{\Frob}{\mathsf{Frob}}
\newcommand{\fp}{\mathfrak{p}}
\newcommand{\fP}{\mathfrak{P}}

\newcommand{\kk}{\mathbb{K}}
\newcommand{\kl}{\mathbb{L}} 

\proof
Assume $p \equiv 1 \md 6$.   Let $\Fp$ be the field with $p$ elements.  We use standard definitions and results from algebraic number theory, as found in \cite{Jan}.  The result essentially follows from the Chebotarev Density Theorem.

Let $\omega$ be a primitive cube root of unity, $\lambda=\sqrt[3]{2}$ be a cube root of $2$ and $\rho=\sqrt[3]{3}$ a cube root of $3$.   Consider the following tower of fields:
$$
M = \Q(\omega, \lambda, \rho) \;\supseteq\; L=\Q(\omega, \lambda)\;\supseteq\;
K=\Q(\omega) \;\supseteq\; \Q.
$$
Let $\O_K$, $\O_L$ denote the rings of integers of $K$ and $L$ respectively.  We may ignore the finitely many ramified primes.  Thus let $p$ be a prime number, sufficiently large that it is unramified in $M$, let $\fp$ be a prime in $K$ extending $p$ and $\fP$ a prime in $L$ extending $\fp$.  Let $\kk=\O_K\!/\fp$ and $\kl=\O_{L}/\fP$ be the residue fields.  We view $\kk$ as embedded in $\kl$ via the map $x+\fp \mapsto x+\fP$.  As $p \equiv 1 \md 6$, $p$ splits in $K$ and $\kk=\O_K\!/\fp \simeq \Fp$.

Since $M$ and $L$ are splitting fields, $M/K$ and $L/K$ are Galois extensions.  The Galois group of $M/K$ is $\textsf{Gal}(M/K)\simeq \Z_3 \times \Z_3$ generated by the maps $\alpha \from \lambda \mapsto \lambda \omega$ and $\beta\from \rho \mapsto \rho \omega$.
The \emph{Frobenius map} of $\kl/\kk$ is the map $x \mapsto x^{|\kl|}$.
The \emph{Frobenius element} $\sigma_{\fp}^L$ is the element of $\textsf{Gal}(L/K)$ inducing the Frobenius map on $\kl/\kk$. (A priori $\sigma_{\fp}^L$ could also depend on the choice of $\fP$ extending $\fp$, but this is not the case since \textsf{Gal}(L/K) is abelian; see \cite[III.2.1]{Jan}.)  Define $\sigma_{\fp}^M \in \textsf{Gal}(M/K)$ analogously.  Then $\sigma_{\fp}^L$ is the restriction of $\sigma_{\fp}^M$ to $L$ by \cite[III.2.3]{Jan}.

By definition of $\kl$, for all sufficiently large $p\equiv 1 \md 6$,
$2\in (\Z_p^*)^3$ if and only if $\kl=\kk$. But $\kl=\kk$ if and only if $\sigma_{\fp}^L$ is the identity map, and it follows that $2\in (\Z_p^*)^3$ if and only if $\sigma_{\fp}^M  \in \langle \beta \rangle$.
Similarly, $3\in (\Z_p^*)^3$ if and only if $\sigma_{\fp}^M  \in \langle \alpha \rangle$ and $6\in (\Z_p^*)^3$ if and only if $\sigma_{\fp}^M  \in \langle \alpha\beta \rangle$.  In summary:
\begin{center}
\begin{tabular}{lcl}
$2,3\notin(\Z_p^*)^3, \ 6\in(\Z_p^*)^3$ &$\iff$ & $\sigma_{\fp}^M\in \{\alpha\beta, \alpha^2\beta^2\}$. \\ %
$2,3,6\notin(\Z_p^*)^3$ &$\iff$ & $\sigma_{\fp}^M\in \{\alpha^2\beta, \alpha\beta^2 \}$. \\ %
\end{tabular}
\end{center}
The Chebotarev Density Theorem \cite[V.10.4]{Jan} implies that for each $\theta \in \textsf{Gal}(M/K)$, the set of primes $\fp$ of $K$  (unramified in $M$) for which $\sigma_{\fp}^M=\theta$ is infinite.  Thus each of the two conditions for $\sigma_{\fp}^M$ displayed above holds infinitely often.\qed

\vspace{0.3cm}

It is possible to describe the primes in Theorem \ref{Thm-p1mod6InfMany} more explicitly.  Given $p\equiv 1 \md 6$, factoring the ideal $p\O_K$ and taking norms, one shows there exist unique $c, d \in \Z$ with $d>0$, $\gcd(c,d)=1$, $c \equiv 2 \md 3$  and $4p=(2c-3d)^2+27d^2$.   Let $t(p)=\left(c \md6,\;  d\md 6\right)$.  There are 9 possible values for $t(p)$: $(2,1)$, $(2,3)$, $(2,5)$, $(5,0)$, $(5,1)$, $(5,2)$, $(5,3)$, $(5,4)$ and $(5,5)$.  The Chebotarev density theorem implies that each of the 9 possible $t(p)$ values occurs ``equally often'' (that is, for a subset of the primes $p \equiv 1 \md 6$ of relative density $1/9$).
Using cubic reciprocity \cite[Ch. 9]{IreRos} one calculates that $2,3\notin(\Z_p^*)^3$ and $6\in(\Z_p^*)^3$ if and only if $t(p)=(2,1)$ or $(5,5)$, while $2,3,6\notin(\Z_p^*)^3$ if and only if $t(p)=(2,5)$ or $(5,1)$.  Each case occurs for $2/9$ of the primes that are $1\md 6$.

The above applications of Lemma \ref{intocircs} have all been with $b=1$. We note however that the conditions of Lemma \ref{intocircs} are never satisfied when $S=\pm\{1,2,3,4\}$ and $b=1$. This is because $2$ is a quadratic residue when $p\equiv 1\md 8$, which means that both $1$ and $4$ are in $H$. The factorisations of $K_p$ into $\cay(\Z_p\,;\pm\{1,2,3,4\})$ in \cite{BrySch} were obtained by applying Lemma \ref{intocircs} with $b=2$ so that $G$ and $H$ have index $2$ and $8$, respectively, in $\Z_p^*$. Another example where Lemma \ref{intocircs} can be applied with $b\neq 1$ is when $p=919$, $S=\pm\{1,2,3\}$, $a=51$ and $b=3$. This yields a factorisation of $K_{919}$ into $\cay(\Z_{919}\,;\pm\{1,2,3\})$. Such a factorisation cannot be obtained by applying Lemma \ref{intocircs} with $b=1$ because $1$, $2$ and $3$ are all cubes in $\Z^*_{919}$.

The following lemma can be used to obtain factorisations of $K_p$, for certain values of $p$, in which some of the factors are isomorphic to $\cay(\Z_p\,;\pm\{1,2,3\})$ and the others are isomorphic to $\cay(\Z_p\,;\pm\{1,2,3,4\})$.

\begin{lemma}\label{mixedfactors}
Let $p$ be prime, let $H$ be the subgroup of $\Z_p^*$ generated by $\{-1,6\}$, and let $d$ be the order of $2H$ in $\Z_p^*/H$. If there exist nonnegative integers $\a$ and $\b$ such that $d=3\a+4\b$, then there is a factorisation of $K_p$ into $\frac{\a(p-1)}{2d}$ copies of $\cay(\Z_p\,;\pm\{1,2,3\})$ and $\frac{\b(p-1)}{2d}$ copies of $\cay(\Z_p\,;\pm\{1,2,3,4\})$.
\end{lemma}

\proof It is sufficient to partition $\Z_p^*$ into $\frac{\a(p-1)}{2d}$ $6$-tuples of the form $\pm\{x,2x,3x\}$ and $\frac{\b(p-1)}{2d}$ $8$-tuples of the form $\pm\{x,2x,3x,4x\}$. Since $d=3\a+4\b$, there is a partition $$\{\{2^{r_i-1}H,2^{r_i}H,2^{r_i+1}H\}:i=1,\ldots,\a\}\cup\{\{2^{r_i-1}H,2^{r_i}H,2^{r_i+1}H,2^{r_i+2}H\}:i=\a+1,\ldots,\a+\b\}$$ of $\{H,2H,\ldots,2^{d-1}H\}$. But $6\in H$ implies $2^{r_i-1}H=3\cdot 2^{r_i}H$ for $i=1,2,\ldots,\a+\b$. Thus, we can rewrite our partition of $\{H,2H,\ldots,2^{d-1}H\}$ as $$\{\{H_i,2H_i,3H_i\}:i=1,\ldots,\a\}\cup\{\{H_i,2H_i,3H_i,4H_i\}:i=\a+1,\ldots,\a+\b\},$$ where $H_i=2^{r_i}H$ for $i=1,\ldots,\a+\b$.

Since $-1\in H$, for $i=1,\ldots,\a$, $H_i\cup 2H_i\cup 3H_i$ can be partitioned into $\frac{|H|}2$ $6$-tuples of the form $\pm\{x,2x,3x\}$, and for $i=\a+1,\ldots,\a+\b$, $H_i\cup 2H_i\cup 3H_i\cup 4H_i$ can be partitioned into $\frac{|H|}2$ $8$-tuples of the form $\pm\{x,2x,3x,4x\}$. If $\mathcal R$ is the set of all $\a\frac{|H|}2$ of these $6$-tuples and $\mathcal S$ is the set of all $\b\frac{|H|}2$ of these $8$-tuples, then $\mathcal R\cup\mathcal S$ is a partition of the subgroup $G=H\cup 2H\cup \cdots\cup 2^{d-1}H$ of $\Z_p^*$. Thus, if $g_1,g_2,\ldots,g_t$ ($t=\frac{p-1}{d|H|}$) represent the cosets of $\Z_p^*/G$, then $$\{g_iR:R\in\mathcal R,i=1,\ldots,t\}\cup\{g_iS:S\in\mathcal S,i=1,\ldots,t\}$$ is a partition of $\Z_p^*$ into $t\a\frac{|H|}2=\frac{\a(p-1)}{2d}$ $6$-tuples of the form $\pm\{x,2x,3x\}$ and $t\b\frac{|H|}2=\frac{\b(p-1)}{2d}$ $8$-tuples of the form $\pm\{x,2x,3x,4x\}$. This is the required partition of $\Z_p^*$.
\qed

\vspace{0.3cm}

Notice that any $6$-factorisation of $K_p$ into $\cay(\Z_p\,;\pm\{1,2,3\})$ given by Lemma \ref{intocircs} can also be obtained via Lemma \ref{mixedfactors}. For if $1,2,3$ represent the three distinct cosets of $G/H$ (where $G=(\Z_p^*)^b$ and $H=(\Z_p^*)^{3b}$, and $p-1=6ab$), then it follows that $\{-1,6\}\subseteq H$ and $2H$ has order $3$ in $G/H$. This means that if $H'$ is the subgroup of $\Z_p^*$ generated by $\{-1,6\}$, then $H'\leq H$ and $3$ divides the order $d$ of $2H'$ in $\Z_p^*/H'$. Thus, we can obtain our $6$-factorisation of $K_p$ into $\cay(\Z_p\,;\pm\{1,2,3\})$ by applying Lemma \ref{mixedfactors} with $\a=\frac d3$ and $\b=0$. Similarly, any $8$-factorisation of $K_p$ into $\cay(\Z_p\,;\pm\{1,2,3,4\})$ given by Lemma \ref{intocircs} can be obtained by applying Lemma \ref{mixedfactors} with $\a=0$ and $\b=\frac d4$.

However, Lemma \ref{mixedfactors} gives us additional factorisations such as the following. When $p=101$ we have $H=\pm\{1,6,14,17,36\}$, and $2H$ has order $d=10$ in $\Z_p^*/H$. Taking $\a=2$ and $\b=1$, we obtain a factorisation of $K_{101}$ into $10$ copies of $\cay(\Z_p\,;\pm\{1,2,3\})$ and $5$ copies of $\cay(\Z_p\,;\pm\{1,2,3,4\})$. Of course, $101$ is neither $1\md 6$ nor $1\md 8$, so there is neither a $6$-factorisation nor an $8$-factorisation of $K_{101}$.

\subsection{Factorising complete graphs of even order}\label{completeminusIintocircs}

In this section we construct factorisations of $K_{2p}-I$ where the factors are all isomorphic to $\cay(\Z_{2p}\,;\pm\{1,2\})$ or all isomorphic to $\cay(\Z_{2p}\,;\pm\{1,2,3,4\})$. We do this by considering $K_{2p}-I$ as a Cayley graph on a dihedral group and partitioning its connection set to generate the factors.
The dihedral group $D_{2p}$ of order $2p$ has elements $r_0,r_1,r_2,\ldots,r_{p-1},s_0,s_1,s_2,\ldots,s_{p-1}$ and satisfies $$r_i\cdot r_j=r_{i+j},\quad r_i\cdot s_j=s_{i+j},\quad s_i\cdot r_j=s_{i-j},\quad s_i\cdot s_j=r_{i-j}$$ where arithmetic of subscripts is carried out modulo $p$.

\begin{lemma}\label{isomorph12}
If $p\geq 3$ is prime, then $$\cay(D_{2p}\,;\{r_{\pm i},s_j,s_{i+j}\})\cong\cay(\Z_{2p}\,;\pm\{1,2\})$$ for all $i\in\Z_p\setminus\{0\}$ and all $j\in\Z_p$.
\end{lemma}

\proof
An isomorphism is given by
$$
\begin{array}{cccccccccccc}
r_0\ &r_i\ &r_{2i}\ &r_{3i}\ &\ldots&r_{(p-1)i}\ &
s_j\ &s_{i+j}\ &s_{2i+j}\ &s_{3i+j}\ &\ldots&s_{(p-1)i+j}\\
\downarrow&\downarrow&\downarrow&\downarrow&\ldots&\downarrow&\downarrow&\downarrow&\downarrow&\downarrow&\ldots&\downarrow\\
0&2&4&6&\ldots&2p-2&2p-1&1&3&5&\ldots&2p-3
\end{array}
$$
\qed

\begin{lemma}\label{isomorph1234}
If $p\geq 5$ is prime, then $$\cay(D_{2p}\,;\{r_{\pm i},r_{\pm 2i},s_j,s_{i+j},s_{2i+j},s_{3i+j}\})\cong\cay(\Z_{2p}\,;\pm\{1,2,3,4\})$$ for all $i\in\Z_p\setminus\{0\}$ and all $j\in\Z_p$.
\end{lemma}

\proof
An isomorphism is given by
$$
\begin{array}{cccccccccccc}
r_0\ &r_i\ &r_{2i}\ &r_{3i}\ &\ldots&r_{(p-1)i}\ &
s_j\ &s_{i+j}\ &s_{2i+j}\ &s_{3i+j}\ &\ldots&s_{(p-1)i+j}\\
\downarrow&\downarrow&\downarrow&\downarrow&\ldots&\downarrow&\downarrow&\downarrow&\downarrow&\downarrow&\ldots&\downarrow\\
0&2&4&6&\ldots&2p-2&2p-3&2p-1&1&3&\ldots&2p-5
\end{array}
$$
\qed

\begin{theorem}\label{K2pminusIinto12}
For each odd prime $p$, there is a factorisation of $K_{2p}-I$ into $\cay(\Z_{2p}\,;\pm\{1,2\})$.
\end{theorem}

\proof
The required factorisation is
$\mathcal F=\{X_i:i\in\Z_p\setminus\{0\}\}$ where $$X_i=\cay(D_{2p}\,;\{r_{\pm 2i},s_i,s_{-i}\})$$ for $i\in\Z_p\setminus\{0\}$. Note that $X_i=X_{-i}$ so $|\mathcal F|=\frac{p-1}2$ as required. Lemma \ref{isomorph12} guarantees that $X_i\cong\cay(\Z_{2p}\,;\pm\{1,2\})$ for each $i\in\Z_p\setminus\{0\}$. Also, $r_0$ is the identity of $D_{2p}$ and each element of $D_{2p}\setminus\{r_0,s_0\}$ occurs in exactly one $X_i$.
Thus, $\mathcal F$ is a factorisation of $\cay(D_{2p}\,;D_{2p}\setminus\{r_0,s_0\})\cong K_{2p}-I$ where the $1$-factor $I$ is $\cay(D_{2p}\,;\{s_0\})$.
\qed

\vspace{0.3cm}

Following work of Davenport \cite[Theorem 5]{Dav} and Weil, a special case of a result due to Moroz
\cite{Mor} yields the following. If $p\equiv 1\md 4$ is prime and $p>8\times 10^6$, then there
exists an integer $x$ such that $x$, $x+1$, $x+2$, $x+3$ represent all four distinct cosets of
$\Z_p^*/(\Z_p^*)^4$. A computer search using {PARI/GP} \cite{pari} verifies in a few minutes that
such an $x$ also exists for all $p< 8\times 10^6$ with $p\equiv 1\md 4$, with the exceptions $p=13$
and $p=17$. Thus, we have the following result.

\begin{lemma}\label{SequencesExist}
If $p\equiv 1\md 4$ is prime with $p\notin\{13,17\}$, then there exists an $x\in\Z_p^*$ such that
$x$, $x+1$, $x+2$ and $x+3$ represent all four distinct cosets of $\Z_p^*/(\Z_p^*)^4$.
\end{lemma}

\begin{theorem}\label{MainTheorem}
If $p\equiv 5\md{8}$ is prime, then
there is a factorisation of $K_{2p}-I$ into $\cay(\Z_{2p}\,;\pm\{1,2,3,4\})$; except that there is no factorisation of $K_{26}-I$ into $\cay(\Z_{2p}\,;\pm\{1,2,3,4\})$.
\end{theorem}

\proof
We first observe that there is no factorisation of $K_{26}-I$ into
$\cay(\Z_{2p}\,;\pm\{1,2,3,4\})$. If such a factorisation exists, then we can assume without loss of generality that the vertex set is
$\Z_{26}$ and that
$\cay(\Z_{26}\,;\pm\{1,2,3,4\})$ is a factor.
But no edge of
$\cay(\Z_{26}\,;\pm\{7\})$ (for example) occurs in a complete subgraph of order $5$ in
$\cay(\Z_{26}\,;\pm\{5,6,7,8,9,10,11,12,13\})$. Since $\cay(\Z_{26}\,;\pm\{1,2,3,4\})$ contains a complete subgraph of order $5$, it follows that
there is no factorisation of $K_{26}-I$ into
$\cay(\Z_{2p}\,;\pm\{1,2,3,4\})$.

Let $p\equiv 5\md 8$ be prime with
$p\neq 13$.
By Lemma \ref{SequencesExist}, there exists an
$x\in\Z_p^*$
such that
$x$, $x+1$, $x+2$ and $x+3$ represent all four distinct cosets of $\Z_p^*/(\Z_p^*)^4$.
By Lemma \ref{isomorph1234},
$$\cay(D_{2p}\,;\{r_{\pm 1},r_{\pm 2},s_x,s_{x+1},s_{x+2},s_{x+3}\})
\cong\cay(\Z_{2p}\,;\pm\{1,2,3,4\}).$$
Now let $H=(\Z_p^*)^4$ act on the subscripts of
the connection set
$\{r_{\pm 1},r_{\pm 2},s_x,s_{x+1},s_{x+2},s_{x+3}\}$
and consider the collection $S_1,S_2,\ldots,S_\frac{p-1}4$
of subsets of $D_{2p}$ thus formed.

We show that
$\{\cay(D_{2p}\,;S_i):i=1,2,\ldots,\frac{p-1}4\}$
is a factorisation of $K_{2p}-I$ into
$\cay(\Z_{2p}\,;\pm\{1,2,3,4\})$.
If $h\in H$, then
$$\cay(D_{2p}\,;\{r_{\pm h},r_{\pm 2h},s_{hx},s_{h(x+1)},s_{h(x+2)},s_{h(x+3)}\})\cong\cay(\Z_{2p}\,;\pm\{1,2,3,4\})$$
by Lemma \ref{isomorph1234} (indeed this is true
for any $h\in\Z_p^*$) so it remains only
to verify that
we have a decomposition of $K_{2p}-I$.
To do this we observe that
$S_1,S_2,\ldots,S_\frac{p-1}4$ partitions
$D_{2p}\setminus\{r_0,s_0\}$
($r_0$ is the identity in $D_{2p}$ and $\cay(D_{2p}\,;\{s_0\})$ is a $1$-factor
in $K_{2p}$).
We have $Hx\cup H(x+1)\cup H(x+2)\cup H(x+3)=\Z_p\setminus\{0\}$.
Also, since $p\equiv 5\md 8$ we have
$-1\in(\Z_p^*)^2$, $-1\notin(\Z_p^*)^4$ and $2\notin(\Z_p^*)^2$
(by the law of quadratic reciprocity).
Thus, $\{\pm h:h\in H\}\cup\{\pm 2h:h\in H\}
=\Z_p\setminus\{0\}$. So
$S_1,S_2,\ldots,S_\frac{p-1}4$ does indeed
partition
$D_{2p}\setminus\{r_0,s_0\}$
and we have the required decomposition.
\qed

\section{2-factorisations of circulant graphs}\label{2FactorisationOfCirculantsSection}

In this section we present various results on $2$-factorisations of circulant graphs, beginning with a couple of known results.
Lemma \ref{2factcir12} was proved independently in \cite{Bry} and \cite{Rod}, and is a special case of a result in \cite{BryMar}. Lemma \ref{2factcirc1234} was proved in \cite{BrySch}.

\begin{lemma}\label{2factcir12}{\rm (\cite{Bry,Rod})}
If $n\geq 5$ and $F$ is any $2$-regular graph of order $n$, then there is a
$2$-factorisation of $\cay(\Z_n\,;\pm\{1,2\})$ into a copy of $F$ and a Hamilton cycle.
\end{lemma}

\begin{lemma}\label{2factcirc1234}{\rm(\cite{BrySch})}
If $n\geq 9$ and $F$ is a $2$-regular graph of order $n$, then there is a $2$-factorisation of $\cay(\Z_n\,;\pm\{1,2,3,4\})$ into $F$ with the definite exceptions of $F=C_4\cup C_5$ and
$F=C_3\cup C_3\cup C_3\cup C_3\cup C_3$,
and the following possible exceptions.
\begin{itemize}
\item [{\rm (1)}] $F=C_3\cup C_3\cup\cdots\cup C_3$ when $n\equiv 3,6\md 9$, $n\geq 21$.
\item [{\rm (2)}] $F=C_4\cup C_4\cup\cdots\cup C_4$ when $n\equiv 4\md 8$, $n\geq 20$.
\item [{\rm (3)}] $F=C_3\cup C_3\cup\cdots\cup C_3\cup C_4$ when $n\equiv 1\md 3$, $n\geq 19$.
\item [{\rm (4)}] $F=C_3\cup C_4\cup C_4\cup\cdots\cup C_4$ when $n\equiv 7\md 8$, $n\geq 23$.
\end{itemize}
\end{lemma}

We now obtain results on $2$-factorisations of $\cay(\Z_n\,;\pm\{1,2,3\})$, but first we need some definitions and notation. For each $m\geq 1$, the graph with vertex set $\{0,1,\ldots,m+2\}$ and edge set $\{\{i,i+1\},\{i+1,i+3\},\{i,i+3\}:i=0,1,\ldots,m-1\}$ is denoted by $J^{1,2,3}_m$. If $F$ is a $2$-regular graph of order $m$, and there exists a decomposition $\{H_1,H_2,H_3\}$ of $J^{1,2,3}_m$ into $F$ such that
\begin{itemize}
\item [(1)] $V(H_1)=\{0,1,\ldots,m+2\}\setminus\{m,m+1,m+2\}$,
\item [(2)] $V(H_2)=\{0,1,\ldots,m+2\}\setminus\{0,2,m+1\}$, and
\item [(3)] $V(H_3)=\{0,1,\ldots,m+2\}\setminus\{0,1,m+2\}$,
\end{itemize}
then we shall write $J^{1,2,3}_m\mapsto F$. Notice that for $i=1,2,3$, the subgraph $H_i$ of $J^{1,2,3}_m$ contains exactly one vertex from each of $\{0,m\}$, $\{1,m+1\}$ and $\{2,m+2\}$.

\begin{lemma}\label{J123givesCay123}
If $n\geq 7$ and $F$ is a $2$-regular graph of order $n$ such that there exists a decomposition $J^{1,2,3}_n\mapsto F$, then there exists a $2$-factorisation of $\cay(\Z_n\,;\pm\{1,2,3\})$ into $F$.
\end{lemma}

\proof
For each $i\in\{0,1,2\}$, identify vertex $i$ of $J^{1,2,3}_n$ with vertex $n+i$. The resulting graph is $\cay(\Z_n\,;\pm\{1,2,3\})$ and the $2$-regular graphs
in the decomposition $J^{1,2,3}_n\mapsto F$ become the required $2$-factors.
\qed

\begin{lemma}\label{J123sumLemma}
If $F$ and $F'$ are vertex-disjoint $2$-regular graphs and there exist decompositions
$J^{1,2,3}_{|V(F)|}\mapsto F$ and $J^{1,2,3}_{|V(F')|}\mapsto F'$, then there exists a
decomposition $J^{1,2,3}_{|V(F)|+|V(F')|}\mapsto F\cup F'$.
\end{lemma}

\proof Let $r$ and $s$ be the respective orders of $F$ and $F'$, let $\{H_1,H_2,H_3\}$ be a
decomposition $J^{1,2,3}_r\mapsto F$ and
let $\{H'_1,H'_2,H'_3\}$ be a decomposition $J^{1,2,3}_s\mapsto F'$. 
Apply the translation $x\mapsto x+r$ to the decomposition $\{H'_1,H'_2,H'_3\}$ to obtain a decomposition
$\{H''_1,H''_2,H''_3\}$ of a copy of $J^{1,2,3}_s$ having vertex set $r,r+1,\ldots,r+s+2$ ($H''_i$ being the translation of $H'_i$ for $i\in\{1,2,3\}$).
It is clear that $\mathcal D=\{H_1\cup H''_1,H_2\cup H''_2,H_3\cup H''_3\}$ is a decomposition
$J^{1,2,3}_{r+s}\mapsto F\cup F'$. Properties (1)-(3) in the definition of $J^{1,2,3}_r\mapsto F$
ensure that $H_i$ and $H''_i$ are vertex-disjoint for $i\in\{1,2,3\}$, and that
\begin{itemize}
\item [(1)] $V(H_1\cup H''_1)=\{0,1,\ldots,r+s+2\}\setminus\{r+s,r+s+1,r+s+2\}$,
\item [(2)] $V(H_2\cup H''_2)=\{0,1,\ldots,r+s+2\}\setminus\{0,2,r+s+1\}$, and
\item [(3)] $V(H_3\cup H''_3)=\{0,1,\ldots,r+s+2\}\setminus\{0,1,r+s+2\}$.
\end{itemize}
\qed

\begin{lemma}\label{123intomcycles}
For each $m\geq 4$, $J^{1,2,3}_m\mapsto C_m$.
\end{lemma}

\proof For $m\in\{4,5,6\}$, $H_1$, $H_2$, $H_3$ are as defined in the following table.

{\small
$$
\begin{array}{|c|c|c|c|}
\hline
m&H_1&H_2&H_3\\
\hline
4
&(0,1,2,3)
&(1,3,6,4)
&(2,4,3,5)\\
\hline
5
&(0,1,2,4,3)
&(1,3,5,7,4)
&(2,3,6,4,5)\\
\hline
6
&(0,1,2,5,4,3)
&(1,3,5,8,6,4)
&(2,4,7,5,6,3)\\
\hline
\end{array}
$$
}

\noindent For $m \geq 7$ and odd
\begin{itemize}
\item $H_1$ contains the edges $\{0,1\}$, $\{1,2\}$, $\{0,3\}$, $\{m-2,m-1\}$ and $\{i,i+2\}$ for $i\in\{2,3,\ldots,m-3\}$,
\item $H_2$ contains the edges $\{1,3\}$, $\{m-2,m\}$, $\{m,m+2\}$, $\{m-1,m+2\}$,  $\{i,i+1\}$ for $i\in\{4,6,\ldots,m-3\}$ and $\{i,i+3\}$ for $i\in\{1,3,\ldots,m-4\}$, and
\item $H_3$ contains the edges $\{2,3\}$, $\{m-2,m+1\}$, $\{m-1,m\}$, $\{m-1,m+1\}$,
$\{i,i+1\}$ for $i\in\{3,5,\ldots,m-4\}$ and $\{i,i+3\}$ for $i\in\{2,4,\ldots,m-3\}$.
\end{itemize}

\noindent For $m \geq 8$ and even
\begin{itemize}
\item $H_1$ contains the edges $\{0,1\}$, $\{1,2\}$, $\{3,4\}$, $\{0,3\}$, $\{2,5\}$, $\{m-2,m-1\}$ and $\{i,i+2\}$ for $i\in\{4,5,\ldots,m-3\}$,
\item $H_2$ contains the edges $\{1,3\}$, $\{1,4\}$, $\{3,5\}$, $\{m-2,m\}$, $\{m,m+2\}$, $\{m-1,m+2\}$, $\{i,i+1\}$ for $i\in\{5,7,\ldots,m-3\}$ and $\{i,i+3\}$ for $i\in\{4,6,\ldots,m-4\}$, and 
\item $H_3$ contains the edges $\{2,4\}$, $\{m-2,m+1\}$, $\{m-1,m\}$, $\{m-1,m+1\}$, $\{i,i+1\}$ for $i\in\{2,4,\ldots,m-4\}$ and $\{i,i+3\}$ for $i\in\{3,5,\ldots,m-3\}$.
\end{itemize}
\qed

\begin{lemma}\label{Additional123Results}
For $m=8$ and for each $m\geq 10$, $J^{1,2,3}_m\mapsto C_3\cup C_{m-3}$.
\end{lemma}

\proof
For $m\in\{8,10,11\}$, $H_1$, $H_2$, $H_3$ are as defined in the following table.

{\small
$$
\begin{array}{|c|c|c|c|}
\hline
m&H_1&H_2&H_3\\
\hline
8
&(4,6,7)\cup(0,1,2,5,3)
&(7,8,10)\cup(1,3,6,5,4)
&(2,3,4)\cup(5,7,9,6,8)\\
\hline
10
&(7,8,9)\cup(0,1,2,4,5,6,3)
&(1,3,4)\cup(5,7,6,9,12,10,8)
&(2,3,5)\cup(4,6,8,11,9,10,7)\\
\hline
11
&(8,9,10)\cup(0,1,2,4,5,7,6,3)
&(1,3,4)\cup(5,6,9,11,13,10,7,8)
&(2,3,5)\cup(4,6,8,11,10,12,9,7)\\
\hline
\end{array}
$$
}

\noindent
For $m\geq 12$ and even
\begin{itemize}
\item $H_1$ consists of the $3$-cycle $(m-3,m-2,m-1)$ and the $(m-3)$-cycle with edges $\{0,1\}$, $\{0,3\}$, $\{1,2\}$, $\{2,4\}$, $\{m-5,m-4\}$, $\{i,i+1\}$ for $i\in\{4,6,\ldots,m-6\}$ and $\{i,i+3\}$ for $i\in\{3,5,\ldots,m-7\}$,
\item $H_2$ consists of the $3$-cycle $(1,3,4)$ and the $(m-3)$-cycle with edges $\{5,7\}$,
    $\{m-5,m-2\}$, $\{m-4,m-3\}$, $\{m-2,m\}$, $\{m,m+2\}$, $\{m-1,m+2\}$, $\{i,i+1\}$ for
    $i\in\{5,7,\ldots,m-7\}$ and $\{i,i+3\}$ for $i\in\{6,8,\ldots,m-4\}$, and
\item $H_3$ consists of the $3$-cycle $(2,3,5)$ and the $(m-3)$-cycle with edges $\{4,6\}$, $\{4,7\}$, $\{m-2,m+1\}$, $\{m-3,m\}$, $\{m-1,m\}$, $\{m-1,m+1\}$ and $\{i,i+2\}$ for $i\in\{6,7,\ldots,m-4\}$.
\end{itemize}

\noindent
For $m\geq 13$ and odd
\begin{itemize}
\item
$H_1$ consists of the $3$-cycle $(m-3,m-2,m-1)$ and the $(m-3)$-cycle with edges $\{0,1\}$, $\{0,3\}$, $\{1,2\}$, $\{2,4\}$, $\{3,6\}$, $\{4,5\}$, $\{5,7\}$, $\{m-5,m-4\}$, $\{i,i+1\}$ for $i\in\{7,9,\ldots,m-6\}$ and $\{i,i+3\}$ for $i\in\{6,8,\ldots,m-7\}$,
\item $H_2$ consists of the $3$-cycle $(1,3,4)$ and the $(m-3)$-cycle with edges $\{5,6\}$,
    $\{m-5,m-2\}$, $\{m-4,m-3\}$, $\{m-2,m\}$, $\{m,m+2\}$, $\{m-1,m+2\}$, $\{i,i+1\}$
    for $i\in\{6,8,\ldots,m-7\}$ and $\{i,i+3\}$ for $i\in\{5,7,\ldots,m-4\}$, and
\item $H_3$ consists of the $3$-cycle $(2,3,5)$ and the $(m-3)$-cycle with edges $\{4,6\}$, $\{4,7\}$, $\{m-2,m+1\}$, $\{m-3,m\}$, $\{m-1,m\}$, $\{m-1,m+1\}$ and $\{i,i+2\}$ for $i\in\{6,7,\ldots,m-4\}$.
\end{itemize}
\qed

\begin{lemma}\label{2factcir123}
Let $n\geq 7$ and let $F$ be a $2$-regular graph of order $n$. If $\nu_3(F)\leq
\nu_5(F)+\sum_{i=7}^{n}\nu_i(F)$ where $\nu_m(F)$ denotes the number of $m$-cycles in $F$, then
there exists a $2$-factorisation of $\cay(\Z_n;\,\pm\{1,2,3\})$ into $F$.
\end{lemma}

\proof
If $n\geq 7$ and $F$ is a $2$-regular graph of order $n$ such that $\nu_3(F)\leq \nu_5(F)+\sum_{i=7}^{n}\nu_i(F)$, then $F$ can be written as a vertex-disjoint union of $2$-regular graphs $G_1,G_2,\ldots,G_t$ where each $G_i$ is isomorphic to either
\begin{itemize}
\item $C_m$ with $m\geq 4$, or
\item $C_3\cup C_{m-3}$ with $m=8$ or $m\geq 10$.
\end{itemize}
By Lemmas \ref{123intomcycles} and \ref{Additional123Results} we have a decomposition $J_{|V(G_i)|}^{1,2,3}\mapsto G_i$ for $i=1,2,\ldots,t$. Applying Lemma \ref{J123sumLemma}
we obtain a decomposition $J^{1,2,3}_n\mapsto F$, and from this we obtain the required
$2$-factorisation of $\cay(\Z_n;\,\pm\{1,2,3\})$ into $F$ by applying Lemma \ref{J123givesCay123}.
\qed

\vspace{0.3cm}

We can obtain an analogue of Lemma \ref{2factcir123} for $\cay(\Z_n\,;\pm\{1,3,4\})$ by using using
similar methods, but we will require $F$ to have girth at least $6$. The graph
with vertex set $\{0,1,\ldots,m+3\}$ and edge set
$\{\{i,i+1\},\{i+1,i+4\},\{i,i+4\}:i=0,1,\ldots,m-1\}$ is denoted by
$J^{1,3,4}_m$. We write $J^{1,3,4}_m\mapsto F$ when there exists a decomposition $\{H_1,H_2,H_3\}$
of $J^{1,3,4}_m$ into a $2$-regular graph $F$ such that
\begin{itemize}
\item [(1)] $V(H_1)=\{0,1,\ldots,m+3\}\setminus\{m,m+1,m+2,m+3\}$,
\item [(2)] $V(H_2)=\{0,1,\ldots,m+3\}\setminus\{0,3,m+1,m+2\}$, and
\item [(3)] $V(H_3)=\{0,1,\ldots,m+3\}\setminus\{0,1,2,m+3\}$.
\end{itemize}
Notice that for $i=1,2,3$, the subgraph $H_i$ of $J^{1,3,4}_m$ contains exactly one vertex from
each of $\{0,m\}$, $\{1,m+1\}$, $\{2,m+2\}$ and $\{3,m+3\}$. It is clear that the proofs of Lemmas
\ref{J123givesCay123} and \ref{J123sumLemma} can be easily modified to give the following two
results.

\begin{lemma}\label{J134givesCay134}
If $n\geq 9$ and $F$ is a $2$-regular graph of order $n$ such that there exists a decomposition
$J^{1,3,4}_n\mapsto F$, then there exists a $2$-factorisation of $\cay(\Z_n\,;\pm\{1,3,4\})$ into
$F$.
\end{lemma}

\begin{lemma}\label{J134sumLemma}
If $F$ and $F'$ are vertex-disjoint $2$-regular graphs and there exist decompositions
$J^{1,3,4}_{|V(F)|}\mapsto F$ and $J^{1,3,4}_{|V(F')|}\mapsto F'$, then there exists a
decomposition $J^{1,3,4}_{|V(F)|+|V(F')|}\mapsto F\cup F'$.
\end{lemma}

Lemmas \ref{J134givesCay134} and \ref{J134sumLemma} allow us to obtain $2$-factorisations of
$\cay(\Z_n\,;\pm\{1,3,4\})$ via the same method we used in the case of $\cay(\Z_n\,;\pm\{1,2,3\})$,
providing we can find appropriate decompositions of $J^{1,3,4}_m$. We now do this.

\begin{lemma}\label{134intomcycles}
For $m=6$, $m=7$ and each $m\geq 9$, $J^{1,3,4}_m\mapsto C_m$.
\end{lemma}

\proof
For $m\in\{6,7,9,10\}$, $H_1$, $H_2$, $H_3$ are as defined in the following table.

{\small
$$
\begin{array}{|c|c|c|c|}
\hline
m&H_1&H_2&H_3\\
\hline
6
&(0,1,5,2,3,4)
&(1,2,6,9,5,4)
&(3,6,5,8,4,7)\\
\hline
7
&(0,1,2,3,6,5,4)
&(1,4,7,10,6,2,5)
&(3,4,8,5,9,6,7)\\
\hline
9
&(0,1,2,3,7,6,5,8,4)
&(1,4,7,8,12,9,6,2,5)
&(3,4,5,9,8,11,7,10,6)\\
\hline
10
&(0,1,2,3,6,9,5,8,7,4)
&(1,4,8,9,13,10,7,6,2,5)
&(3,4,5,6,10,9,12,8,11,7)\\
\hline
\end{array}
$$
}

\noindent
For $m \geq 11$ and odd
\begin{itemize}
\item $H_1$ contains the edges $\{0,1\}$, $\{0,4\}$, $\{1,2\}$, $\{2,3\}$, $\{3,7\}$, $\{5,6\}$, $\{m-3,m-2\}$, $\{m-5,m-1\}$, $\{m-4,m-1\}$ and $\{i,i+4\}$ for $i\in\{4,5,\ldots,m-6\}$,
\item $H_2$ contains the edges $\{1,4\}$, $\{1,5\}$, $\{2,5\}$, $\{2,6\}$, $\{4,7\}$, $\{m,m+3\}$, $\{m-1,m+3\}$, $\{m-2,m-1\}$, $\{m-3,m\}$, $\{i,i+1\}$ for $i\in\{7,9,\ldots,m-4\}$ and $\{i,i+3\}$ for $i\in\{6,8,\ldots,m-5\}$, and
\item $H_3$ contains the edges $\{3,4\}$, $\{3,6\}$, $\{4,5\}$, $\{m-1,m\}$, $\{m-2,m+1\}$, $\{m-1,m+2\}$, $\{m-4,m\}$, $\{m-3,m+1\}$, $\{m-2,m+2\}$,
$\{i,i+1\}$ for $i\in\{6,8,\ldots,m-5\}$ and $\{i,i+3\}$ for $i\in\{5,7,\ldots,m-6\}$.
\end{itemize}
For $m \geq 12$ and even
\begin{itemize}
\item $H_1$ contains the edges $\{0,1\}$, $\{0,4\}$, $\{1,2\}$, $\{2,3\}$, $\{3,6\}$, $\{4,7\}$, $\{5,6\}$, $\{5,9\}$, $\{m-5,m-2\}$, $\{m-4,m-3\}$, $\{m-4,m-1\}$, $\{m-2,m-1\}$, $\{i,i+1\}$ for $i\in\{7,9,\ldots,m-7\}$ and $\{i,i+3\}$ for $i\in\{8,10,\ldots,m-6\}$,
\item $H_2$ contains the edges $\{1,4\}$, $\{1,5\}$, $\{2,5\}$, $\{2,6\}$, $\{4,8\}$, $\{m-6,m-2\}$, $\{m-5,m-4\}$, $\{m-5,m-1\}$, $\{m-3,m-2\}$, $\{m-3,m\}$, $\{m-1,m+3\}$, $\{m,m+3\}$, $\{i,i+1\}$ for $i\in\{6,8,\ldots,m-8\}$ and $\{i,i+3\}$ for $i\in\{7,9,\ldots,m-7\}$, and
\item $H_3$ contains the edges $\{3,4\}$, $\{3,7\}$, $\{4,5\}$, $\{5,8\}$, $\{6,9\}$, $\{m-6,m-5\}$, $\{m-4,m\}$, $\{m-3,m+1\}$, $\{m-2,m+1\}$, $\{m-2,m+2\}$, $\{m-1,m\}$, $\{m-1,m+2\}$ and $\{i,i+4\}$ for $i\in\{6,7,\ldots,m-7\}$.
\end{itemize}
\qed

\begin{lemma}\label{Additional134Results}
For each $m\geq 14$, $J^{1,3,4}_m\mapsto C_8\cup C_{m-8}$.
\end{lemma}

\proof
For $m\in\{14,15,16,17\}$, $H_1$, $H_2$, $H_3$ are as defined in the following table.

{\small
$$
\begin{array}{|c|l|}
\hline
m&\\
\hline
14
&H_1=(0,1,2,3,7,8,5,4)\cup(6,9,13,12,11,10)\\
&H_2=(8,11,14,17,13,10,9,12)\cup(1,4,7,6,2,5)\\
&H_3=(7,10,14,13,16,12,15,11)\cup(3,4,8,9,5,6)\\
\hline
15
&H_1=(0,1,2,3,6,5,8,4)\cup(7,10,14,13,9,12,11)\\
&H_2=(1,4,7,8,9,6,2,5)\cup(10,11,14,18,15,12,13)\\
&H_3=(8,11,15,14,17,13,16,12)\cup(3,4,5,9,10,6,7)\\
\hline
16
&H_1=(0,1,5,6,2,3,7,4)\cup(8,9,10,11,15,14,13,12)\\
&H_2=(1,2,5,9,6,7,8,4)\cup(10,13,16,19,15,12,11,14)\\
&H_3=(3,4,5,8,11,7,10,6)\cup(9,12,16,15,18,14,17,13)\\
\hline
17
& H_1=(0,1,2,3,7,6,5,4)\cup(8,9,13,16,12,15,14,10,11)\\
& H_2=(1,4,8,12,9,6,2,5)\cup(7,10,13,14,17,20,16,15,11)\\
& H_3=(3,4,7,8,5,9,10,6)\cup(11,12,13,17,16,19,15,18,14)\\
\hline
\end{array}
$$
}

\noindent
For $m\geq 18$ and even
\begin{itemize}
\item $H_1$ consists of the $8$-cycle $(0,1,5,6,2,3,7,4)$ and the $(m-8)$-cycle with edges $\{8,9\}$, $\{9,10\}$, $\{10,11\}$, $\{8,12\}$, $\{m-5,m-1\}$, $\{m-4,m-3\}$, $\{m-3,m-2\}$,
$\{m-2,m-1\}$
 $\{i,i+1\}$ for $i\in\{12,14,\ldots,m-6\}$ and $\{i,i+3\}$ for $i\in\{11,13,\ldots,m-7\}$,
\item $H_2$ consists of the $8$-cycle $(1,2,5,9,6,7,8,4)$ and the $(m-8)$-cycle with edges $\{10,13\}$, $\{11,12\}$, $\{m-6,m-2\}$, $\{m-5,m-2\}$, $\{m-4,m-1\}$, $\{m-3,m\}$, $\{m-1,m+3\}$, $\{m,m+3\}$ and $\{i,i+4\}$ for $i\in\{10,11,\ldots,m-7\}$, and
\item $H_3$ consists of the $8$-cycle $(3,4,5,8,11,7,10,6)$ and the $(m-8)$-cycle with edges $\{9,12\}$, $\{9,13\}$, $\{m-4,m\}$, $\{m-3,m+1\}$, $\{m-2,m+1\}$, $\{m-2,m+2\}$, $\{m-1,m\}$, $\{m-1,m+2\}$, $\{i,i+1\}$ for $i\in\{13,15,\ldots,m-5\}$ and $\{i,i+3\}$ for $i\in\{12,14,\ldots,m-6\}$.
\end{itemize}

\noindent
For $m\geq 19$ and odd
\begin{itemize}
\item
$H_1$ consists of the $8$-cycle $(0,1,2,3,7,6,5,4)$ and the $(m-8)$-cycle with edges $\{8,9\}$, $\{8,11\}$, $\{9,13\}$, $\{10,11\}$, $\{10,14\}$, $\{12,15\}$, $\{12,16\}$, $\{m-4,m-1\}$, $\{m-3,m-2\}$ and $\{i,i+4\}$ for $i\in\{13,14,\ldots,m-5\}$,
\item $H_2$ consists of the $8$-cycle $(1,4,8,12,9,6,2,5)$ and the $(m-8)$-cycle with edges $\{7,10\}$, $\{7,11\}$, $\{10,13\}$, $\{11,15\}$, $\{m-4,m-3\}$, $\{m-3,m\}$, $\{m-2,m-1\}$, $\{m-1,m+3\}$, $\{m,m+3\}$, $\{i,i+1\}$ for $i\in\{13,15,\ldots,m-6\}$ and $\{i,i+3\}$ for $i\in\{14,16,\ldots,m-5\}$, and
\item $H_3$ consists of the $8$-cycle $(3,4,7,8,5,9,10,6)$ and the $(m-8)$-cycle with edges $\{11,12\}$, $\{11,14\}$, $\{12,13\}$, $\{m-4,m\}$, $\{m-3,m+1\}$, $\{m-2,m+1\}$, $\{m-2,m+2\}$, $\{m-1,m\}$, $\{m-1,m+2\}$, $\{i,i+1\}$ for $i\in\{14,16,\ldots,m-5\}$ and $\{i,i+3\}$ for $i\in\{13,15,\ldots,m-6\}$.
\end{itemize}
\qed

\begin{lemma}\label{134into888}
$J^{1,3,4}_{24}\mapsto C_8\cup C_8\cup C_8$.
\end{lemma}

\proof
Take
\begin{align*}
H_1&=(0,1,2,3,6,5,8,4)\cup (7,10,9,12,13,14,15,11)\cup(16,17,18,19,23,22,21,20),\\
H_2&=(1,4,7,8,9,6,2,5)\cup (10,11,12,15,16,13,17,14)\cup (18,21,24,27,23,20,19,22), \mbox{ and} \\
H_3&=(3,4,5,9,13,10,6,7)\cup (8,11,14,18,15,19,16,12)\cup (17,20,24,23,26,22,25,21).
\end{align*}
\qed

\vspace{0.3cm}

The following result is an analogue of Lemma \ref{2factcir123} for $2$-factorisations of $\cay(\Z_n\,;\pm\{1,3,4\})$.

\begin{lemma}\label{2factcir134}
If $n\geq 9$ and $F$ is a $2$-regular graph of order $n$ with girth at least $6$, 
then there exists a $2$-factorisation of $\cay(\Z_n\,;\pm\{1,3,4\})$ into $F$.
\end{lemma}

\proof
If $n\geq 9$ and $F$ is a $2$-regular graph of order $n$ with girth at least $6$, then $F$ can be written as a vertex-disjoint union of $2$-regular graphs $G_1,G_2,\ldots,G_t$ where each $G_i$ is isomorphic to either
\begin{itemize}
\item $C_m$ with $m=6,7$ or  $m\geq 9$, 
\item $C_8\cup C_{m-8}$ with $m\geq 14$, or
\item $C_8\cup C_8\cup C_8$.
\end{itemize}
By Lemmas \ref{134intomcycles}, \ref{Additional134Results} and \ref{134into888} we have a decomposition $J_{|V(G_i)|}^{1,3,4}\mapsto G_i$ for $i=1,2,\ldots,t$. Applying Lemma \ref{J134sumLemma}
we obtain a decomposition $J^{1,3,4}_n\mapsto F$, and from this we obtain the required
$2$-factorisation of $\cay(\Z_n;\,\pm\{1,3,4\})$ into $F$ by applying Lemma \ref{J134givesCay134}.
\qed

\section{$2$-factorisations and the Oberwolfach Problem}

In this section we use results from the preceding sections to obtain results on the Oberwolfach Problem (and an additional result on $2$-factorisations of $K_n-I$ into a number of specified $2$-factors and Hamilton cycles). We will also use the following corollary of Lemma \ref{2factcirc1234} which was proved in \cite{BrySch}.

\begin{lemma}\label{CircDecompImpliesAll}{\rm (\cite{BrySch})}
If there exists a factorisation of $K_n$ or of $K_n-I$ into $\cay(\Z_n\,;\pm\{1,2,3,4\})$, then $\op(F)$ has a solution for each $2$-regular graph $F$ of order $n$, with the exception that there is no solution to $\op(C_4\cup C_5)$.
\end{lemma}

\begin{theorem}\label{mainThmp=5mod8}
If $p\equiv 5\md 8$ is prime, then $\op(F)$ has a solution for every $2$-regular graph $F$ of order $2p$.
\end{theorem}

\proof
The case $p=13$ is covered in \cite{DezFraHuaMesRos}. For $p\neq 13$, Theorem \ref{MainTheorem} gives us a factorisation of $K_{2p}-I$ into $\cay(\Z_{2p}\,;\pm\{1,2,3,4\})$ and the result then follows by Lemma \ref{CircDecompImpliesAll}.
\qed

\begin{theorem}\label{combinedFactorsThm}
Let $\P$ be the set of primes given by $p\in\P$ if and only if $p \geq 7$ and neither $4$ nor $32$
is in the subgroup of $\Z_p^*$ generated by $\{-1,6\}$. Then $\P$ is infinite and if $p\in\P$, then
$\op(F)$ has a solution for every $2$-regular graph $F$ of order $p$ satisfying
$\nu_3(F)\leq\nu_5(F)+\sum_{i=7}^{n}\nu_i(F)$ where $\nu_m(F)$ denotes the number of $m$-cycles in
$F$.
\end{theorem}

\proof
Let $p$ be prime such $p\equiv 1\md 6$, $2,3\notin(\Z_p^*)^3$ 
and $6\in(\Z_p^*)^3$. Theorem \ref{Thm-p1mod6InfMany} says that there are infinitely many such $p$. 
We shall show that $p\in\P$, which shows that $\P$ is also infinite.
We have  $-1\in(\Z_p^*)^3$, and this together with the fact that 
$6\in(\Z_p^*)^3$ implies that the subgroup of $\Z_p^*$ generated by $\{-1,6\}$ is a subgroup of $(\Z_p^*)^3$.
Since it follows from $2\notin(\Z_p^*)^3$ that $4,32\notin(\Z_p^*)^3$, neither $4$ nor $32$
is in the subgroup of $\Z_p^*$ generated by $\{-1,6\}$. That is, $p\in\P$.  

Now let $p$ be an arbitrary element of $\P$ and let $G$ be the subgroup of $\Z_p^*$ generated by $\{-1,6\}$. The condition that
neither $4$ nor $32$ is in $G$ implies that the order $d$ of $2G$ in $\Z_p^*/G$ is neither $1$, $2$
nor $5$, and so there exist non-negative integers $\a$ and $\b$ such that $d=3\a+4\b$. Thus, by
Lemma \ref{mixedfactors} there is a factorisation of $K_p$ in which each factor is either
$\cay(\Z_p\,;\pm\{1,2,3\})$ or $\cay(\Z_p\,;\pm\{1,2,3,4\})$.

Let $F$ be a $2$-regular graph of order $p$ satisfying
$\nu_3(F)\leq\nu_5(F)+\sum_{i=7}^{n}\nu_i(F)$. Lemma \ref{2factcir123} gives us a $2$-factorisation
of $\cay(\Z_p\,;\pm\{1,2,3\})$ into $F$, and Lemma \ref{2factcirc1234} gives us a $2$-factorisation
of $\cay(\Z_p\,;\pm\{1,2,3,4\})$ (the facts that $p$ is prime and that
$\nu_3(F)\leq\nu_5(F)+\sum_{i=7}^{n}\nu_i(F)$ imply that $F$ is not amongst the possible exceptions
listed in Lemma \ref{2factcirc1234}). The result follows. \qed

\begin{theorem}\label{134result}
Let $\P$ be the set of primes such that $p\in\P$ if and only if $p\equiv 1\md 6$ and
$2,3,6\notin(\Z_p^*)^3$. Then $\P$ is infinite and if $p\in\P$, then $\op(F)$ has a solution for every $2$-regular
graph $F$ of order $p$ with girth at least $6$.
\end{theorem}

\proof By Theorem \ref{Thm-p1mod6InfMany}, $\P$ is infinite. If $p \in \mathcal{P}$, then Theorem
\ref{Kpinto134} gives us a factorisation of $K_p$ into $\cay(\Z_p\,;\pm\{1,3,4\})$, and the result 
then follows by applying Lemma \ref{2factcir134} to each factor ($7\notin\P$ so Lemma \ref{2factcir134} can indeed be applied). \qed

\vspace{0.3cm}

For each odd prime $p$, the following theorem states there is a $2$-factorisation of $K_{2p}-I$ into $\frac{p-1}2$ prescribed $2$-factors and $\frac{p-1}2$ Hamilton cycles.

\begin{theorem}
If $p$ is an odd prime and $G_1,G_2,\ldots,G_{\frac{p-1}{2}}$ are $2$-regular graphs of order $2p$, then there is a $2$-factorisation $\{F_1,F_2,\ldots,F_{p-1}\}$ of $K_{2p}-I$ such that $F_i\cong G_i$ for $i=1,2,\ldots,\frac{p-1}{2}$ and $F_i$ is a Hamilton cycle for $i=\frac{p+1}2,\frac{p+3}2,\ldots,p-1$.
\end{theorem}

\proof
By Theorem \ref{K2pminusIinto12} there is a factorisation of $K_{2p}-I$ into $\cay(\Z_p\,;\pm\{1,2\})$. By Lemma \ref{2factcir12}, each copy of $\cay(\Z_p\,;\pm\{1,2\})$ can be factored into any specified $2$-regular graph of order $2p$ and a Hamilton cycle. The result follows.
\qed

\section{Isomorphic 2-factorisations of complete multigraphs}

The complete multigraph of order $n$ and multiplicity $s$ is denoted by $sK_n$. It has $s$ distinct edges joining each pair of distinct vertices.

\begin{lemma}\label{lambdafoldintocirc}
If $p$ is an odd prime and $S=\pm\{d_1,d_2,\ldots,d_s\}\subseteq\Z_p^*$, then there exists a $2s$-factorisation of $sK_p$ into $\cay(\Z_p\,;S)$.
\end{lemma}

\proof
The required factorisation is given by
$\{\cay(\Z_p\,;\omega^iS):i=0,1,\ldots,\frac{p-3}2\}$
where $\omega$ is primitive in $\Z_p$ and $\omega^iS=\{\omega^is:s\in S\}$.
\qed

\begin{theorem}\label{3Kpthm}
If $p$ is an odd prime and $F$ is any $2$-regular graph of order $p$ satisfying
$\nu_3(F)\leq\nu_5(F)+\sum_{i=7}^{n}\nu_i(F)$, where $\nu_m(F)$ denotes the number of $m$-cycles in
$F$, then there exists a $2$-factorisation of $3K_p$ into $F$.
\end{theorem}

\proof
The cases $p=3$ and $p=5$ are trivial so assume $p\geq 7$.
By Lemma \ref{lambdafoldintocirc} there exists a $6$-factorisation of $3K_p$ into $\cay(\Z_p\,;\pm\{1,2,3\})$, and by Lemma \ref{2factcir123} each such $6$-factor has a $2$-factorisation into $F$.
\qed

\begin{theorem}\label{4Kpthm}
If $p$ is an odd prime and $F$ is any $2$-regular graph of order $p$,
then there exists a $2$-factorisation of $4K_p$ into $F$.
\end{theorem}

\proof The cases $p=3$ and $p=5$ are trivial. Since solutions to $\op(C_7)$ and $\op(C_3 \cup C_4)$
exist, the case $p=7$ can be dealt with by taking four copies of these $2$-factorisations of $K_7$.
So we may assume $p\geq 11$. By Lemma \ref{lambdafoldintocirc} there exists an $8$-factorisation of
$4K_p$ into $\cay(\Z_p\,;\pm\{1,2,3,4\})$, and by Lemma \ref{2factcirc1234} each such $8$-factor
has a $2$-factorisation into $F$; except in the case where $F$ is one of the listed exceptions or
possible exceptions in Lemma  \ref{2factcirc1234}. These are easily dealt with as follows. Since
$p$ is prime the only relevant exceptions are $F=C_3\cup C_3\cup\cdots\cup C_3\cup C_4$ where the
number of copies of $C_3$ is at least $5$, and $F=C_3\cup C_4\cup C_4\cup\cdots\cup C_4$ where the
number of copies of $C_4$ is odd and at least $5$. However, it is known that for each such $F$,
there is a $2$-factorisation of $K_p$ into $F$; the former case is covered in \cite{DejFraMenRos},
and the latter case is covered in \cite{Koh1}. Thus, by taking four copies of these
$2$-factorisations of $K_p$, we obtain the required $2$-factorisations of $4K_p$. \qed

\begin{theorem}\label{lambdageq7thm}
Let $p$ be an odd prime and let $F$ be a $2$-regular graph of order $p$. If $\lambda\equiv 0\md 4$,
then there exists a $2$-factorisation of $\lambda K_p$ into $F$. Moreover, if $F$ satisfies
$\nu_3(F)\leq\nu_5(F)+\sum_{i=7}^{n}\nu_i(F)$, where $\nu_m(F)$ denotes the number of $m$-cycles in
$F$, then the result also holds for $\lambda=3$ and for all $\lambda\geq 6$.
\end{theorem}

\proof
For the given values of $\lambda$, it is trivial to factorise $\lambda K_p$ such that each factor is
either $3K_p$ or $4K_p$, and with each factor being $4K_p$ when $\lambda\equiv 0\md 4$. Thus,
the result follows by Theorems \ref{3Kpthm} and \ref{4Kpthm}.
\qed

\vspace{0.3cm} \noindent{\bf Acknowledgement.}
The authors acknowledge the support of the Australian Research Council via grants DE120100040, DP0770400, DP120100790, DP120103067 and DP130102987.

\end{document}